\documentclass[12pt,a4paper]{amsart}

\usepackage[T1]{fontenc}
\usepackage{lmodern}

\newtheorem{theorem}{Theorem}
\newtheorem{lemma}[theorem]{Lemma}
\newtheorem{e-proposition}[theorem]{Proposition}

\newtheorem{e-definition}[theorem]{Definition\rm}

\DeclareMathAlphabet{\doba}{U}{msb}{m}{n}         

\gdef\mR{\doba{R}}

\def\e{\varepsilon}

\def\cF{\mathcal{F}}
\def\cL{\mathcal{L}}

\def\E{\mathcal{E}}   
\def\Ric{{\mathop{\rm Ric}}}

\def\di{\mathrm{d}}
\def\tn{\tilde\nabla}
\def\bD{\bar\Delta}
\def\bn{\overline\nabla}

\let\<\langle 
\let\>\rangle

\newcommand{\setR}{\mathbb{R}}

\def\og{\leavevmode\raise.3ex\hbox{$\scriptscriptstyle\langle\!\langle$~}}
\def\fg{\leavevmode\raise.3ex\hbox{~$\!\scriptscriptstyle\,\rangle\!\rangle$}}

\begin{document}
\centerline{}

\title{A construction of conformal-harmonic maps}

\author[Biquard]{Olivier Biquard}
\author[Madani]{Farid Madani}

\address{UPMC Universit\'e Paris 6 and \'Ecole Normale Sup\'erieure}
\address{Fakult\"at f\"ur  Mathematik, Universit\"at Regensburg, Germany}

\begin{abstract}

Conformal harmonic maps from a 4-dimensional conformal manifold to a
Riemannian manifold are maps satisfying a certain conformally
invariant fourth order equation.  We prove a general existence result
for conformal harmonic maps, analogous to the Eells-Sampson theorem
for harmonic maps. The proof uses a geometric flow and relies on
results of Gursky-Viaclovsky and Lamm.
\end{abstract}

\maketitle

Let $(M,g)$ and $(N,h)$ be two compact Riemannian manifolds of
dimension 4 and $n$ respectively.  Denote by $R^M$, $\Ric^M$ and $S^M$
the Riemann, Ricci and scalar curvatures associated to $(M,g)$
respectively.  The tension field $\tau (u)\in u^*TN$ of a map $u\in C^2(M,N)$
is defined by:
$$\tau(u)=-\sum \tn_{e_i} \di u(e_i),$$
where $\{e_i\}$ is an orthonormal frame of $TM$ and $\tilde\nabla$, $\bn$ are
the Riemannian connections on $T^*M\otimes u^*TN$ and $u^*TN$
respectively. So $u$ is a harmonic map if and only if $\tau(u)=0$. 

The following conformal energy functional for maps $u:M\to N$ was
introduced by B\'erard \cite{Berard}:
\begin{equation}\label{energie}
\E(u)=\int_M\big(|\tau(u)|^2+\frac{2}{3}S^M|\di u|^2-2\Ric^M(\di u,\di u)\big)\di v_g
\end{equation}
It is conformally invariant with respect to conformal changes of $g$
(but not of $h$). The map $u$ is said to be a conformal-harmonic map
(in short, C-harmonic map), if it is a critical point of $\E$. Namely,
it is a solution of the following equation:
\begin{equation}\label{ech}
\cL(u):=\bD\tau(u)+R^N(\di u(e_i), \tau(u))\di u(e_i)+\bn^*\{(\frac{2}{3}S^M-2\Ric^M)\di u\}=0
\end{equation}
where $\bD=\bn^*\bn$ is the rough Laplacian, $\bn^*$ is
the $L^2$-adjoint of $\bn$. This C-harmonic
equation \eqref{ech} differs from the biharmonic equation by low
order terms which make the C-harmonic equation conformally
invariant.

The Yamabe number $\mu(M,[g])$ and the total $Q$-curvature $\kappa(M,[g])$
are defined by
\begin{equation}
 \mu(M,[g])= \inf_{g'\in [g]} \frac{\int_M S_{g'}\di v_{g'}}{(\int_M\di v_{g'})^{1/2}},\qquad \kappa(M,[g])=\frac{1}{12}\int_M \big(S_g^2-3|\Ric_g|^2\big)\di v_{g}.
\end{equation}
Both are conformal invariants of $[g]$.
The aim of this Note is to prove the following:
\begin{theorem}\label{exist}
  Let $(M^4,[g])$ and $(N^n,h)$ be compact manifolds equipped with a
  conformal metric $[g]$ and a Riemannian metric $h$ respectively.
  Assume that the curvature of $(N,h)$ is non positive, $\mu(M,[g])>0$
  and $\kappa(M,[g])+\frac{1}{6}\mu^2(M, [g])> 0$. Then each homotopy class
  in $C^\infty(M,N)$ can be represented by a C-harmonic map.
\end{theorem} 
This statement is analogous to the classical theorem of Eells and
Sampson on the existence of harmonic maps, but no uniqueness is
proved. Remark that C-harmonic maps are usually not harmonic for any
metric in the conformal class $[g]$, so the theorem actually
constructs new applications which depend only on the conformal
geometry of $[g]$.

The rest of the paper is devoted to the proof of Theorem \ref{exist}.
The geometric hypothesis is explained by the following crucial
proposition, which is adapted from the theorem of Gursky and
Viaclovsky \cite{GV} giving conditions for the Paneitz operator to be
positive (this corresponds to the case $N=\setR$). The conditions on
$(M,[g])$ are the same as in \cite{GV}, but we add a non positive
curvature assumption on $(N,h)$:
\begin{e-proposition}\label{Gtheo}
  Let $(M,g)$ and $(N,h)$ be compact Riemannian manifolds of dimension
  4 and $n$ respectively. Assume that the curvature of $(N,h)$ is
  non positive, $\mu(M,[g])>0$ and $\kappa(M,[g])+\frac{1}{6}\mu^2(M, [g])>
  0$. Then there exists a positive constant $c$ such that $\E(u)\geq
  c\|\di u\|_2^2$.
\end{e-proposition}
It is important to note that the RHS of the inequality is not conformally
invariant, so the constant $c$ depends on the metric $g$
itself. Nevertheless, since the LHS is conformally invariant, it is
sufficient to prove the inequality just for one special metric in the
conformal class.

\noindent\emph{Proof.}
  We adapt \cite[\S~6]{GV}. For any map $u\in C^\infty(M,N)$, we have the
  Bochner-Weitzenb\"ock formula:
  \begin{equation}
  \|\tau(u)\|_2^2=\|\tn\di u\|_2^2\\-\int_M \bigl(\sum_{i,j}\langle R^N_{\di u(e_i),\di
    u(e_j)}\di u(e_j),\di u(e_i)\rangle-\Ric^M(\di u,\di u) \bigr) \di v_g.\label{Boch}
  \end{equation}
  We combine \eqref{energie} and $\frac{4}{3}$ times \eqref{Boch} and
  use the fact that the curvature of $N$ is non positive, to get
  \begin{equation}
  \E(u)\geq \frac{4}{3}\big(\|\tn\di u\|_2^2-\frac{1}{4}\|\tau(u)\|_2^2\big)+\frac{2}{3}\int_M \big(S^M|\di u|^2-\Ric^M(\di u,\di u)\big)\di v_g.
  \end{equation}
  Again by \cite{GV}, under the hypothesis, there exists a metric $g$
  in the conformal class such that the scalar curvature and the second
  symmetric function of the eigenvalues of $\Ric$ are positive, which
  implies $ \Ric^M<S^M g$. Combining with the fact that $|\tn\di
  u|_2^2-\frac{1}{4}|\tau(u)|_2^2=|\tn_0\di u|^2$, where $\tn_0\di u$ is
  the trace free part of $\tn\di u$, we deduce that there exists $c>0$
  such that $\E(u)\geq c\|\di u\|_2^2$. \qed

In order to prove Theorem \ref{exist}, we study the gradient flow
of the energy:
\begin{equation}\label{gfe}
\partial_tu=-\cL(u), \quad 
u(\cdot, 0)=u_0.
\end{equation}
Equation \eqref{gfe} is a fourth order strongly parabolic
equation. So, there exists $T>0$ such that equation \eqref{gfe} has
a solution $u\in C^\infty(M\times [0,T),N)$. We have the following estimates:
\begin{lemma}\label{estimate}
  Under the assumptions of Theorem \ref{exist} and if $u\in C^\infty(M\times
  [0,T),N)$ is a solution of \eqref{gfe}, then we have for all
  $t\in[0,T)$:
\begin{align}
\E(u(\cdot,t))+ 2 \int_0^t\|\partial_tu\|^2_2\di t&=\E(u_0),\label{fond}\\
\|\di u(\cdot,t)\|_2^2+2\int_0^t\|\bn\tau(u)\|_2^2\di t&\leq \| \di u_0\|_2^2+ct\label{2fond},\\
\int_M|\tn\di u|^2\di v_g&\leq c\label{3fond},\\
\int_M|\di u|^4\di v_g&\leq c\label{4fond}.
\end{align}
\end{lemma}

\noindent\emph{Proof.}
Let $u$ be a solution of \eqref{gfe}. We have 
$\frac{\di}{\di t}\E(u)=\di\E (\partial_t u)=2\langle \cL(u),\partial_tu\rangle=-2\langle \partial_tu,\partial_tu\rangle$.
Equality \eqref{fond} holds after an integration over time.

Using Proposition \ref{Gtheo} and \eqref{fond}, we obtain
$ \|\di u\|_2\leq c$.
Using again \eqref{fond}, we conclude that 
$ \|\tau (u)\|_2\leq c$.
Using \eqref{Boch} and the fact that $N$ has non positive curvature, we have 
$\|\tn\di u\|_2^2\leq \|\tau (u)\|_2^2+c  \|\di u\|_2^2$.
Hence we obtain \eqref{3fond}.

Now observe that
$$ \cL(u) = \bD\tau(u)+R^N(\di u(e_i), \tau(u))\di u(e_i)+(\nabla R^M)\odot\di u+R^M\odot \tilde\nabla \di u, $$
where the $\odot$'s are fixed bilinear operations. Since $N$ has
nonnegative curvature, using (\ref{3fond}) and the bounds on $\|\di
u\|_2$ and $\|\tau (u)\|_2$, it follows that
$ \langle \cL(u) ,\tau(u)\rangle \geq \|\bn\tau(u)\|_2^2 - c $
for some constant $c$. Therefore,
\begin{equation*}
\frac12\frac{\di}{\di t}\|\di u(\cdot,t)\|_2^2 =\langle \partial_tu,\tau(u)\rangle=-\langle \cL(u) ,\tau(u)\rangle \leq -\|\bn\tau(u)\|_2^2+c .
\end{equation*} 
Hence \eqref{2fond} holds after an integration over time. Inequality
\eqref{4fond} is obtained from the bound on $\|\di u\|_2$ by using
\eqref{3fond} and Sobolev embedding.\qed

Given the estimates of the lemma, the arguments for proving the long
time existence of the flow and its convergence are very similar to
that developed by Lamm \cite{Lam} in its study of the flow for
biharmonic maps. Below we briefly explain how to use its proof to finish the
proof of Theorem \ref{exist}.

Biharmonic maps are solutions of the equation
\begin{equation}\label{biharmonic}
\bD\tau(u)+R^N(\di u(e_i), \tau(u))\di u(e_i)=0
\end{equation}
which are critical points of the functional 
\begin{gather}
\cF(u)=\int_M|\tau(u)|^2\di v_g .
\end{gather}
If we compare the functionals $\E$, $\cF$, we
note that there exists $c>0$ such that $$\cF(u)-c\|\di u\|_2^2\leq\E(u)\leq
\cF(u)+c\|\di u\|_2^2.$$ Since $\| \di u\|_2$ is estimated in Lemma
\ref{estimate}, all the estimates obtained by Lamm for the solutions
of the biharmonic map gradient flow hold for the $C-$harmonic
ones. Therefore, we give a series of results where the proofs are
given in \cite{Lam} (itself building on deep regularity results of
Chang-Yang and Struwe). To state them, we need to introduce
a smooth nonnegative cut-off function $\eta$ on $M$, such that 
\begin{equation}
  \begin{split}
    \eta&=1 \text{ on } B_R(P_0),\\ \eta&=0 \text{ on } M-B_{2R}(P_0),\\
    \|\nabla^i\eta\|_\infty&\leq \frac{c}{R^{2i}} \text{ for } i=1 \text{ and } 2,
  \end{split}
\end{equation}
where $B_R(P_0)$ is a geodesic ball with centre $P_0\in M$ and  radius $R<\min(1,\frac{1}{4}\mathrm{inj}_M(P_0))$.
Let $E$ be the (conformally invariant) total energy defined by
\begin{equation}
E(u)=\E(u)+\biggl(\int_M|\di u|^4\di v_g\biggr)^\frac{1}{2}.
\end{equation}
The local energy $E(u;B_{R}(P_0))$ is defined by an analogous formula as $E$, where the integration is taken over the ball $B_R(P_0)$.
\begin{lemma}\label{abcd}
 Let $u\in C^\infty(M\times [0,T),N)$ be a solution of \eqref{gfe}. There exists $\e_0=\e_0(u_0,M,N)$ such that if 
$\sup_{0\leq t<T} E(u(\cdot,t);B_{2R}(P_0))<\e_0$, then for all $t<T$ and $k=2,3$,
$$ \int_0^t\int_M\eta^4|\tn^k\di u|^2\di v_g\di t\leq c+\frac{ct}{R^4}.$$
\end{lemma}

\begin{lemma}\label{no sing}
If $v:\mR^4\longrightarrow N$ is a weak C-harmonic map 
with $E(v)<+\infty$, then $E(v)=0$.
\end{lemma}

\begin{lemma}\label{rere}
 Let $u\in C^\infty(M\times [0,T),N)$ be a solution of \eqref{gfe}. There exists $\e_0=\e_0(u_0,M,N)$ such that if 
$\sup_{(P,t)\in M\times[0,T)} E(u(\cdot,t);B_{2R}(P))<\e_0$, then there exists $\delta\in (0,\min(T,cR^4))$ such that the H\"older norms of $u$ and all derivatives of $u$ are uniformly bounded on $[t-\frac{\delta}{4},t]$ for all $t\in[\frac{\delta}{2},T)$ by constants which depend only on $u_0$, $M,\;N,\;R$ and the order of derivatives of $u$.
\end{lemma}

Given these three lemmas, we now finish the proof of Theorem
\ref{exist}. This still follows Lamm, but it may be useful for the
reader to follow how one quickly deduces the result from the previous
technical lemmas. As mentioned above, there exists $T>0$ and $u\in
C^\infty(M\times [0,T),N)$ solution of \eqref{gfe}.  Two cases can occur:
\begin{enumerate}
\item for all $R>0$, there exists $P_0\in M$ such that
$$\limsup_{t\to T^-}E(u(\cdot,t);B_{R}(P_0))\geq\e_0;$$
\item there exists $\e_0$ and $R>0$ such that
$$\sup_{(P,t)\in M\times [0,T)}E(u(\cdot,t);B_{R}(P))<\e_0.$$
\end{enumerate}
The first case is when we have concentration. 
Then we claim that there exists a finite number of singularities. We
construct a sequence of maps $\tilde u_m:\mR^4\longrightarrow N$, by
$  \tilde u_m (x)=u(P_m+R_mx,t_m+R^4_ms_m)$,
where $(R_m), (t_m)$ and $(P_m)$ are sequences of positive real
numbers and points of $M$ which converge to 0, $T$ and $P_0\in M$
respectively. There exists a sequence $(s_m)$ with values in $[t_0,
0]$, such that $\tilde u_m\to \tilde u$ weakly in
$W^{4,2}_{loc}(\mR^4,N)$ and strongly in $W^{3,2}_{loc}(\mR^4,N)$. It
implies that $\tilde u$ is a weak C-harmonic map, with finite and
positive total energy $E(\tilde u)$. From the regularity Lemma \ref{no sing}, we get
a contradiction.

So we must be in the second case: Lemma \ref{rere} then implies
that the flow can be continued smoothly beyond $T$.  Therefore we
conclude that the heat flow doesn't have any singularities and exists
for all times.  From Lemma \ref{abcd} and \eqref{fond}, we conclude that
there exists $\tau>0$ such that
\begin{equation*}\int^{t+\tau}_t\int_M|\nabla^4u|^2\di v_g\di t\leq c
\quad\text{ and }\quad\lim_{t\to+\infty}\int^{t+\tau}_t\int_M|\partial_t u|^2\di v_g\di t=0.
\end{equation*}
There exists $t_m\to+\infty$ such that $u_m:=u(\cdot,t_m)$ converges to $u_\infty$
weakly in $W^{4,2}(M,N)$ and $\partial_t u_m$ converges to 0 in
$L^2(M,N)$. By Lemma \ref{rere}, we have that $u_m$ is uniformly
bounded in $C^k$ for all $k\in\mathbb N$. Taking a limit in \eqref{gfe},
we deduce that $u_\infty$ is a smooth C-harmonic map.





\end{document}